\documentclass[12pt]{amsart}
\usepackage{mathrsfs}
\usepackage{cases}
\usepackage{epic}
\usepackage{amsfonts}
\usepackage{amsmath}
\usepackage{amssymb, upgreek}
\usepackage{bm}
\usepackage{latexsym,todonotes}
\usepackage{pdflscape}
\usepackage[all]{xypic}
\usepackage[all]{xy}
\usepackage{color}
\usepackage{colordvi}
\usepackage{multicol}
\usepackage[normalem]{ulem}
\usepackage[linktocpage=true]{hyperref}
\hypersetup{colorlinks,linkcolor=blue,urlcolor=cyan,citecolor=blue}

\begin{document}
	\input xy
	\xyoption{all}
	
	\numberwithin{equation}{section}
	\allowdisplaybreaks
	\renewcommand{\mod}{\operatorname{mod}\nolimits}
	\newcommand{\proj}{\operatorname{proj.}\nolimits}
	\newcommand{\rad}{\operatorname{rad}\nolimits}
	\newcommand{\soc}{\operatorname{soc}\nolimits}
	\newcommand{\ind}{\operatorname{inj.dim}\nolimits}
	\newcommand{\id}{\operatorname{id}\nolimits}
	\newcommand{\Mod}{\operatorname{Mod}\nolimits}
	\newcommand{\R}{\operatorname{R}\nolimits}
	\newcommand{\End}{\operatorname{End}\nolimits}
	\newcommand{\ob}{\operatorname{Ob}\nolimits}
	\newcommand{\Ht}{\operatorname{Ht}\nolimits}
	\newcommand{\cone}{\operatorname{cone}\nolimits}
	\newcommand{\rep}{\operatorname{rep}\nolimits}
	\newcommand{\Ext}{\operatorname{Ext}\nolimits}
	\newcommand{\Tor}{\operatorname{Tor}\nolimits}
	\newcommand{\Hom}{\operatorname{Hom}\nolimits}
	\newcommand{\Pic}{\operatorname{Pic}\nolimits}
	\newcommand{\aut}{\operatorname{Aut}\nolimits}
	\newcommand{\Fac}{\operatorname{Fac}\nolimits}
	\newcommand{\Div}{\operatorname{Div}\nolimits}
	\newcommand{\rank}{\operatorname{rank}\nolimits}
	\newcommand{\Len}{\operatorname{Length}\nolimits}
	\newcommand{\RHom}{\operatorname{RHom}\nolimits}
	\renewcommand{\deg}{\operatorname{deg}\nolimits}
	\renewcommand{\Im}{\operatorname{Im}\nolimits}
	\newcommand{\Ker}{\operatorname{ker}\nolimits}
	\newcommand{\Iso}{\operatorname{Iso}\nolimits}
	\newcommand{\Coh}{\operatorname{coh}\nolimits}
	\newcommand{\Qcoh}{\operatorname{Qch}\nolimits}
	\newcommand{\inj}{\operatorname{inj.dim}\nolimits}
	\newcommand{\dimv}{\operatorname{\underline{\dim}}\nolimits}
	\newcommand{\res}{\operatorname{res}\nolimits}
	\def \TT{\mathbf{T}}
	\def \U{\mathbf U}
	\def \tUB{{}^{\texttt{B}}\tU}
	\def \UB{{}^{\texttt{B}}\U}
	\newcommand{\indim}{\operatorname{inj.dim}\nolimits}
	\newcommand{\pdim}{\operatorname{proj.dim}\nolimits}
	\def \dbl{\rm dbl}
	\def \sqq{\mathbf v}
	\newcommand{\Aut}{\operatorname{Aut}\nolimits}
	
	\def \bm{\mathbf m}

	\newcommand{\Cp}{\operatorname{Cp}\nolimits}
	\newcommand{\coker}{\operatorname{Coker}\nolimits}
	\renewcommand{\dim}{\operatorname{dim}\nolimits}
	\renewcommand{\ker}{\operatorname{Ker}\nolimits}
	\renewcommand{\div}{\operatorname{div}\nolimits}
	\newcommand{\Ab}{{\operatorname{Ab}\nolimits}}
	\newcommand{\Cone}{{\operatorname{Cone}\nolimits}}
	\renewcommand{\Vec}{{\operatorname{Vec}\nolimits}}
	\newcommand{\pd}{\operatorname{proj.dim}\nolimits}
	\newcommand{\sdim}{\operatorname{sdim}\nolimits}
	\newcommand{\add}{\operatorname{add}\nolimits}
	\newcommand{\pr}{\operatorname{pr}\nolimits}
	\newcommand{\oR}{\operatorname{R}\nolimits}
	\newcommand{\oL}{\operatorname{L}\nolimits}
	
	\newcommand{\tU}{\operatorname{\widetilde{\bf U}}\nolimits}
	
	\newcommand{\Perf}{{\mathfrak Perf}}
	\newcommand{\cc}{{\mathcal C}}
	\newcommand{\ce}{{\mathcal E}}
	\newcommand{\cs}{{\mathcal S}}
	\newcommand{\cf}{{\mathcal F}}
	\newcommand{\cx}{{\mathcal X}}
	\newcommand{\cy}{{\mathcal Y}}
	\newcommand{\cl}{{\mathcal L}}
	\newcommand{\ct}{{\mathcal T}}
	\newcommand{\cu}{{\mathcal U}}
	\newcommand{\cm}{{\mathcal M}}
	\newcommand{\cv}{{\mathcal V}}
	\newcommand{\ch}{{\mathcal H}}
	\newcommand{\ca}{{\mathcal A}}
	\newcommand{\mcr}{{\mathcal R}}
	\newcommand{\cb}{{\mathcal B}}
	\newcommand{\ci}{{\mathcal I}}
	\newcommand{\cj}{{\mathcal J}}
	\newcommand{\cp}{{\mathcal P}}
	\newcommand{\cg}{{\mathcal G}}
	\newcommand{\cw}{{\mathcal W}}
	\newcommand{\co}{{\mathcal O}}
	\newcommand{\cd}{{\mathcal D}}
	\newcommand{\ck}{{\mathcal K}}
	\newcommand{\calr}{{\mathcal R}}
	
	\def \fg{{\mathfrak g}}
	\newcommand{\ol}{\overline}
	\newcommand{\ul}{\underline}
	\newcommand{\cz}{{\mathcal Z}}
	\newcommand{\st}{[1]}
	\newcommand{\ow}{\widetilde}
	\renewcommand{\P}{\mathbf{P}}
	\newcommand{\pic}{\operatorname{Pic}\nolimits}
	\newcommand{\Spec}{\operatorname{Spec}\nolimits}
	\newtheorem{theorem}{Theorem}[section]
	\newtheorem{acknowledgement}[theorem]{Acknowledgement}
	\newtheorem{algorithm}[theorem]{Algorithm}
	\newtheorem{axiom}[theorem]{Axiom}
	\newtheorem{case}[theorem]{Case}
	\newtheorem{claim}[theorem]{Claim}
	\newtheorem{conclusion}[theorem]{Conclusion}
	\newtheorem{condition}[theorem]{Condition}
	\newtheorem{conjecture}[theorem]{Conjecture}
	\newtheorem{construction}[theorem]{Construction}
	\newtheorem{corollary}[theorem]{Corollary}
	\newtheorem{criterion}[theorem]{Criterion}
	\newtheorem{propdef}[theorem]{Definition-Proposition}
	
	\newtheorem{definition}[theorem]{Definition}
	\newtheorem{example}[theorem]{Example}
	\newtheorem{exercise}[theorem]{Exercise}
	\newtheorem{lemma}[theorem]{Lemma}
	\newtheorem{notation}[theorem]{Notation}
	\newtheorem{problem}[theorem]{Problem}
	\newtheorem{proposition}[theorem]{Proposition}
	\newtheorem{solution}[theorem]{Solution}
	\newtheorem{summary}[theorem]{Summary}
	\newtheorem*{thm}{Theorem}
	\newcommand{\qbinom}[2]{\begin{bmatrix} #1\\#2 \end{bmatrix} }

	\theoremstyle{remark}
	\newtheorem{remark}[theorem]{Remark}
	
	\def \bfk{\mathbf k}
	\def \bp{{\mathbf p}}
	\def \bA{{\mathbf A}}
	\def \bL{{\mathbf L}}
	\def \bF{{\mathbf F}}
	\def \bS{{\mathbf S}}
	\def \bC{{\mathbf C}}
	\def \bD{{\mathbf D}}
	\def \Ire{\I^{\rm re}}
	\def \Iim{\I^{\rm im}}
	\def \Iiso{\I^{\rm iso}}
	
	\def \bs{\mathbf s}
	\def \I{\mathbb I}
	
	\def \Z{{\Bbb Z}}
	\def \F{{\Bbb F}}
	\def \C{{\Bbb C}}
	\def \N{{\Bbb N}}
	\def \Q{{\Bbb Q}}
	\def \G{{\Bbb G}}
	\def \X{{\Bbb X}}
	\def \P{{\Bbb P}}
	\def \K{{\Bbb K}}
	\def \E{{\Bbb E}}
	\def \A{{\Bbb A}}
	\def \BH{{\Bbb H}}
	\def \T{{\Bbb T}}
	\newcommand{\bluetext}[1]{\textcolor{blue}{#1}}
	\newcommand{\redtext}[1]{\textcolor{red}{#1}}
	\newcommand{\red}[1]{\redtext{ #1}}
	\newcommand{\blue}[1]{\bluetext{ #1}}
	\def \tMH{{\cs\cd\widetilde{\ch}}}
	
	\title[Quantum Borcherds-Bozec algebras via semi-derived Ringel-Hall algebras II]{Quantum Borcherds-Bozec algebras via semi-derived Ringel-Hall algebras II: braid group actions}
	
	\author{Ji Lin}
	\address{Department of Mathematics and Statistics, Fuyang Normal University, Fuyang 236037, P.R.China}
	\email{jlin@fynu.edu.cn}

	\author[Ming Lu]{Ming Lu}
	\address{Department of Mathematics, Sichuan University, Chengdu 610064, P.R.China}
	\email{luming@scu.edu.cn}

	\author[Shiquan Ruan]{Shiquan Ruan}
	\address{ School of Mathematical Sciences,
		Xiamen University, Xiamen 361005, P.R.China}
	\email{sqruan@xmu.edu.cn}

	\dedicatory{Dedicated to Professor Liangang Peng on the occasion of his 65th birthday}

	\subjclass[2010]{Primary 17B37, 
		16E60, 18E35.}
	\keywords{Quantum Borcherds-Bozec algebras, Semi-derived Ringel-Hall algebras, Braid group actions, Quivers with loops}
	
	\begin{abstract}
		Based on the realization of quantum Borcherds-Bozec algebra $\tU$ and quantum generalized Kac-Moody algebra $\tUB$ via semi-derived Ringel-Hall algebra of a quiver with loops, we deduce the braid group actions of $\tU$ introduced by Fan and Tong recently and establish braid group actions for $\tUB$ by applying the BGP reflection functors to semi-derived Ringel-Hall algebras.
	\end{abstract}
	
	\maketitle
	\section{Introduction}
	
	\subsection{}
	
	Quantum groups (i.e., quantized enveloping algebras) are introduced by Drinfeld and Jimbo in 1980s.
	Ringel \cite{Rin90} in 1990 constructed a Hall algebra associated to a Dynkin quiver $Q$ over a finite field $\mathbb F_q$, and identified its generic version with half part of  a quantum group $\U^+ =\U^+_v(\fg)$; see Green \cite{Gr} for an extension to quivers {\em without loops}. Ringel's construction has led to a geometric construction of $\U^+$ by Lusztig, who in addition constructed its canonical basis \cite{L90}. These constructions can be regarded as the earliest examples of categorifications of halves of quantum groups.
	
	Bridgeland \cite{Br13} in 2013 succeeded in using a Hall algebra of complexes to realize the quantum group $\U$. Actually Bridgleland's construction naturally produces the Drinfeld double $\widetilde \U$, a variant of $\U$ with the Cartan subalgebra doubled (with generators $K_i, K_i'$, for $i \in \I$). A  reduced version, which is the quotient of $\widetilde \U$ by the ideal generated by the central elements $K_i K_i'-1$, is then identified with $\U$. 
	
	Bridgeland's version of Hall algebras has found further generalizations and improvements which allow more flexibilities. More recently, motivated by the works of Bridgeland and Gorsky, the second author and Peng \cite{LP16} formulated the {\em semi-derived Ringel-Hall algebras} starting from hereditary abelian categories (including module categories over arbitrary quivers); compare with the semi-derived Hall algebras  defined in \cite{Gor18} for Frobenius categories.  
	
	\subsection{}
	Associated to a Borcherds-Cartan matrix, the {\em generalized Kac-Moody algebra} is introduced by Borcherds in 1988. Its quantization, called {\em quantum generalized Kac-Moody algebra}, is defined by Kang \cite{K95}.
	The Ringel-Hall algebra realization of the positive half $\UB^+$ of a quantum generalized Kac-Moody algebra is given in \cite{KS06}, where they also constructed the canonical basis.
	
	\subsection{}
	The {\em Borcherds-Bozec algebra} is a further generalization of the Kac-Moody algebra, which is also  associated to a Borcherds-Cartan matrix, and
	generated by infinitely many generators. 
	The {\em quantum Borcherds-Bozec algebra} is introduced by Bozec in research of the perverse sheaves of a quiver with loops \cite{B15,B16}, which can be viewed as a further generalization of the quantum generalized Kac-Moody algebra. It is generated by $K_i^{\pm1},e_{il},f_{il}$ ($(i,l)\in\I^\infty$), where $\I$ is a finite or countably infinite index set, and $\I^\infty$ is defined in \S\ref{subsec: QBB}. Bozec also constructed the canonical basis \cite{B16}. The Ringel-Hall algebra realization of the positive half $\U^+$ of a quantum Borcherds-Bozec algebra is given in \cite{K18} by considering a quiver with loops.

	Thanks to Bozec, there exists a set of primitive
	generators $\mathfrak{a}_{il}, \mathfrak{b}_{il}$ ($(i, l)\in\I^\infty$) with better properties and simpler commutation relations.
	Using these generators, Bozec then developed
	the crystal basis theory for quantum Borcherds-Bozec algebras \cite{B16}. The
	global bases for the quantum Borcherds-Bozec algebras are constructed
	in \cite{FKKT20}.
	
	Recently, Fan and Tong \cite{FT21} generalized Lusztig's braid group actions of quantum groups to quantum Borcherds-Bozec algebras by using the primitive
	generators $\mathfrak{a}_{il}, \mathfrak{b}_{il}$ ($(i, l)\in\I^\infty$).

	\subsection{}
	In \cite{Lu22}, the second author used the semi-derived Ringel-Hall algebras of quivers with loops to realize the whole quantum Borcherds-Bozec algebras and quantum generalized Kac-Moody algebras, a generalization of Brigeland's construction. Similar to the quantum groups, the construction produces the Drinfeld double $\tU$ and $\tUB$, with their Cartan subalgebras doubled. The quotients of $\tU$ and $\tUB$ by the ideal generated by the central elements $K_i K_i'-1$ give the  quantum Borcherds-Bozec algebra $\U$ and quantum generalized Kac-Moody algebra $\UB$.

	\subsection{}
	In this paper, we shall use the isomorphisms of semi-derived Ringel-Hall algebras induced by the BGP reflection functors of quivers to deduce the braid group actions of quantum Borcherds-Bozec algebras introduced by Fan and Tong in \cite{FT21}. As the set of generators $e_{el},f_{il}$ ($(i,l)\in\I^\infty$) are interpreted clearly in Hall algebras, we first give the formulas of braid group actions on these generators, and then prove these formulas coincide with
	the ones given in \cite{FT21}. 
	
	We also use the isomorphisms of semi-derived Ringel-Hall algebras induced by the BGP reflection functors of quivers to give the braid group actions $T'_{i,e}$, $T''_{i,e}$ ($e=\pm1$) of quantum generalized Kac-Moody algebras. Indeed, we prove that $T'_{i,e}$, $T''_{i,e}$
	are algebra isomorphisms and they satisfy the defining relations of the braid groups associated to Borcherds-Cartan matrices.
	
	\subsection{}
	
	The paper is organized as follows. In Section \ref{sec:SDRHA}, we review the basic materials on Hall algebras and semi-derived Ringel-Hall algebras. Section \ref{sec:Drinfeld} is devoted to deducing the braid group actions of quantum Borcherds-Bozec algerbras introduced in \cite{FT21} by using the BGP reflection functors. In Section \ref{sec:Quantum generalized}, we construct the braid group actions of quantum generalized Kac-Moody algebras.
	
	\vspace{2mm}
	\noindent{\bf Acknowledgments.}
	JL is partially supported by the National Natural Science Foundation of China (No. 12001107) and the Natural Science Foundation of the Anhui Higher Education Institutions of China (No. KJ2021A0661, KJ2021A0658). ML is partially supported by the National Natural Science Foundation of China (No. 12171333). SR is partially supported by the National Natural Science Foundation of China (No. 12271448) and the Fundamental Research Funds for Central Universities of China (No. 20720220043). We thank the anonymous referee for very helpful suggestions and comments.

	\section{Semi-derived Ringel-Hall algebras}
	\label{sec:SDRHA}

	Denote by $\Z_2:=\Z/2\Z$.
	In this section, we shall review categories of $\Z_2$-graded complexes, Ringel-Hall algebras, and semi-derived Ringel-Hall algebras for arbitrary hereditary abelian categories over a finite field $\bfk=\F_q$.
	
	\subsection{Categories of $\Z_2$-graded complexes}
	We assume that $\ca$  is an exact category. 
	
	Let $\cc_{\Z_2}(\ca)$ be the exact category of $\Z_2$-graded complexes over $\ca$. Namely, an object $M$ of this category is a diagram with objects and morphisms in $\ca$:
	$$\xymatrix{ M^0 \ar@<0.5ex>[r]^{d^0}& M^1 \ar@<0.5ex>[l]^{d^1}  },\quad d^1d^0=d^0d^1=0.$$
	All indices of components of $\Z_2$-graded objects will be understood modulo $2$.
	A morphism $s=(s^0,s^1):M\rightarrow N$ is a diagram
	\[\xymatrix{  M^0 \ar@<0.5ex>[r]^{d^0} \ar[d]^{s^0}& M^1 \ar@<0.5ex>[l]^{d^1} \ar[d]^{s^1} \\
		N^0 \ar@<0.5ex>[r]^{e^0}& N^1 \ar@<0.5ex>[l]^{e^1} }  \]
	with $s^{i+1}d^i=e^is^i$.
	
	The shift functor on complexes is an involution
	\begin{align}
		\label{inv}
		\xymatrix{\cc_{\Z_2}(\ca) \ar[r]^{*} & \cc_{\Z_2}(\ca)\ar[l],}
	\end{align}
	which shifts the grading and changes the sign of the differential as follows:
$$\xymatrix{ (\;M^0 \ar@<0.5ex>[r]^{d^0}& M^1\;) \ar@<0.5ex>[l]^{d^1} \ar[r]^{*} & (\;M^1 \ar@<0.5ex>[r]^{-d^1} \ar[l]& M^0\;) \ar@<0.5ex>[l]^{-d^0}  }.$$
	
	Let $\bfk$ be a field.
	In the following, we always assume that $\ca$ is a hereditary $\bfk$-linear abelian  category which is essentially small with finite-dimensional homomorphism and extension spaces.
	
	Similar to ordinary complexes, we can define the $i$-th homology group for $M$, denoted by $H^i(M)$, for any $i\in\Z_2$. A complex is called acyclic if its homology group is zero. The subcategory formed by all acyclic complexes is denoted by $\cc_{\Z_2,ac}(\ca)$.
	
	For any object $X\in\ca$, we define
	\begin{align}
		\label{stalks}
		\begin{split}
			K_X:=&(\xymatrix{ X \ar@<0.5ex>[r]^{1}& X \ar@<0.5ex>[l]^{0}  }),\qquad \,\, K_X^*:=(\xymatrix{ X \ar@<0.5ex>[r]^{0}& X \ar@<0.5ex>[l]^{1}  }),
			\\
			C_X:=&(\xymatrix{ 0 \ar@<0.5ex>[r]& X \ar@<0.5ex>[l]  }),\qquad \quad C_X^*:=(\xymatrix{ X\ar@<0.5ex>[r]& 0 \ar@<0.5ex>[l]  })
		\end{split}
	\end{align}
	in $\cc_{\Z_2}(\ca)$. Note that $K_X,K_X^*$ are acyclic complexes.
	

	\begin{lemma}[\text{\cite[Proposition 2.3]{LP16}}]
		\label{proposition extension 2 zero}
		Let $\ca$ be a hereditary abelian category.
		For any $K\in\cc_{\Z_2,ac}(\ca)$,  we have
		$$\pdim_{\cc_{\Z_2}(\ca)} (K)\leq 1\text{ and }\indim_{\cc_{\Z_2}(\ca)} (K)\leq 1.$$
	\end{lemma}

	Denote by $\Iso(\cc_{\Z_2}(\ca))$ the set of isomorphism classes $[M]$ of
	$\cc_{\Z_2}(\ca)$.
	
	By Proposition \ref{proposition extension 2 zero}, for any $K,M\in \cc_{\Z_2}(\ca)$ with $K$ acyclic, define
	\begin{align*}
		\langle K,M\rangle=\dim\Hom_{\cc_{\Z_2}({\ca})}(K,M)-\dim\Ext^1_{\cc_{\Z_2}({\ca})}(K,M),
		\\
		\langle M,K\rangle =\dim\Hom_{\cc_{\Z_2}({\ca})}(M,K)-\dim\Ext^1_{\cc_{\Z_2}({\ca})}(M,K).
	\end{align*}
	We call them the Euler forms. They descend to bilinear forms on the Grothendieck groups $K_0(\cc_{\Z_2,ac}(\ca))$
	and $K_0(\cc_{\Z_2}(\ca))$, again denoted by the same symbol $\langle\cdot,\cdot\rangle$. 

	Let $K_0(\ca)$ be the Grothendieck group of $\ca$.
	For any $A\in\ca$, we denote by $\widehat{A}$ the corresponding element in  $K_0(\ca)$.
	We also use $\langle \cdot,\cdot\rangle$ to denote the Euler form of $\ca$, i.e.,
	\begin{align*}
	\langle A,B\rangle:=	\langle \widehat{A}, \widehat{B}\rangle=\dim\Hom_\ca(A,B)-\dim\Ext^1_\ca(A,B), \text{ for any }A,B\in\ca.
	\end{align*}
	Let $(\cdot,\cdot)$ be the symmetrized Euler form of $\ca$, i.e., $(A,B)= \langle A, B\rangle+\langle B, A\rangle$.
	
	
	\begin{proposition}[\text{\cite[Proposition 2.4, Corollary 2.5]{LP16}}]
		\label{lema euler form}
		For any $A ,B \in\ca$, we have the following.
		\begin{align}
			&\langle C_A,K_B\rangle=\langle C_A^*,K_B^*\rangle=\langle \widehat{A}, \widehat{B}\rangle,\quad
			\langle K_B,C_A^*\rangle=\langle K_B^*,C_A\rangle=\langle \widehat{B},\widehat{A}\rangle;
			\\
			&\langle K_B,C_A\rangle=\langle C_A^*,K_B\rangle=\langle C_A,K_B^*\rangle=\langle K_B^*,C_A^*\rangle=0;
			\\
			&\langle K_{A}, K_{B}\rangle= \langle K_{A}^*, K_{B}^*\rangle=\langle K_{A}, K_{B}^*\rangle=\langle K_{A}^*, K_{B}\rangle=\langle \widehat{A},\widehat{B}\rangle.
		\end{align}
	\end{proposition}

	Let $\ca\coprod\ca$ be the product of two copies of $\ca$. Then there is a forgetful functor $\res: \cc_{\Z_2}(\ca)\rightarrow \ca\coprod \ca$, which maps $M=(\xymatrix{ M^0 \ar@<0.5ex>[r]^{d^0}& M^1 \ar@<0.5ex>[l]^{d^1}  })$ to $(M^0, M^{1})$.

	
	
	\subsection{Ringel-Hall algebras}
	
	Let $\ca$ be an essentially small abelian category, linear over the finite field $\bfk=\F_q$.

	Given objects $X,Y,Z\in\ca$, define $\Ext^1_\ca(X,Z)_Y\subseteq \Ext^1_\ca(X,Z)$ to be the subset parameterising extensions with the middle term  isomorphic to $Y$. We define the Ringel-Hall algebra (also called Hall algebra) $\ch(\ca)$ to be the $\Q$-vector space whose basis is formed by the isomorphism classes $[X]$ of objects $X$ of $\ca$, with the multiplication
	defined by
	\begin{align}
		\label{eq:mult}
		[X]\diamond [Z]=\sum_{[Y]\in \Iso(\ca)}\frac{|\Ext_\ca^1(X,Z)_Y|}{|\Hom_\ca(X,Z)|}[Y].
	\end{align}
	It is well known that
	the algebra $\ch(\ca)$ is associative and unital. The unit is given by $[0]$, where $0$ is the zero object of $\ca$; see \cite{Rin90,Br13}. 

	For any three objects $X,Y,Z$, let
	\begin{align}
		\label{eq:Fxyz}
		F_{XZ}^Y:= \big |\{L\subseteq Y \mid L \cong Z,  Y/L\cong X\} \big |.
	\end{align}
	The Riedtmann-Peng formula states that
	\[
	F_{XZ}^Y= \frac{|\Ext^1(X,Z)_Y|}{|\Hom(X,Z)|} \cdot \frac{|\Aut(Y)|}{|\Aut(X)| |\Aut(Z)|}.
	\]
	
	For any object $X$, 
	let
	\begin{align*}
		[\![X]\!]:=\frac{[X]}{|\Aut(X)|}.
	\end{align*}
	Then the Hall multiplication \eqref{eq:mult} can be reformulated to be
	\begin{align}
		[\![X]\!]\diamond [\![Z]\!]=\sum_{[\![Y]\!]}F_{X,Z}^Y[\![Y]\!],
	\end{align}
	which is the version of Hall multiplication used in \cite{Rin90}.

	\subsection{Semi-derived Ringel-Hall algebras}
	
	Let 
	$$\sqq=\sqrt{q}.$$
	Let $\widetilde{\ch}(\cc_{\Z_2}(\ca))$ be the twisted Ringel-Hall algebra of $\cc_{\Z_2}(\ca)$ over $\Q(\sqq)$, that is, $\widetilde{\ch}(\cc_{\Z_2}(\ca))$ has a basis formed by the isomorphism classes $[M]$ of objects $M$ of $\cc_{\Z_2}(\ca)$, with the product given by
	\begin{align}\label{multiplication formula}
		[L]* [M]=&\sqq^{\langle \res L,\res M\rangle}[L]\diamond [M]
		\\\notag
		=&\sqq^{\langle \res L,\res M\rangle}\sum_{[X]\in \Iso(\cc_{\Z_2}(\ca))}\frac{|\Ext^1_{\cc_{\Z_2}(\ca)}(L,M)_X|}{|\Hom_{\cc_{\Z_2}(\ca)}(L,M)|}[X].
	\end{align}
	
	Let $I_{\Z_2}$ be the two-sided ideal of $\widetilde{\ch}(\cc_{\Z_2}(\ca))$ generated by all differences $[L]-[K\oplus M]$ if there is a short exact sequence
	\begin{equation}
		\label{eq:ideal}
		0 \longrightarrow K \longrightarrow L \longrightarrow M \longrightarrow 0
	\end{equation}
	with $K$ acyclic. 

	
	Let $\widetilde{\ch}(\cc_{\Z_2}(\ca))/I_{\Z_2}$ be the quotient algebra. 
	We also denote by $*$ the induced multiplication in $\widetilde{\ch}(\cc_{\Z_2}(\ca))/I_{\Z_2}$. In the following, we shall use the same symbols both in $\widetilde{\ch}(\cc_{\Z_2}(\ca))$ and $\widetilde{\ch}(\cc_{\Z_2}(\ca))/I_{\Z_2}$.
	Let
	\begin{align}
		\label{multiset}
		S_{\Z_2}:=\{a[K] \in \widetilde{\ch}(\cc_{\Z_2}(\ca))/I_{\Z_2} \mid  a \in \Q(\sqq)^\times, K\in\cc_{\Z_2,ac}(\ca)\},
	\end{align}
	which is a multiplicatively closed subset with the identity $[0]\in S_{\Z_2}$.

	\begin{propdef}[\cite{LP16}]
		\label{proposition localizaition of Ringel-Hall algebra}
		The multiplicatively closed subset $S_{\Z_2}$ is a right Ore, right reversible subset of $\widetilde{\ch}(\cc_{\Z_2}(\ca))/I_{\Z_2}$. Equivalently, the right localization of $\widetilde{\ch}(\cc_{\Z_2}(\ca))/I_{\Z_2}$ with respect to $S_{\Z_2}$ exists, which is called the twisted semi-derived Ringel-Hall algebra of $\ca$, and denoted by $\tMH(\ca)$.
	\end{propdef}

	
	For any $\alpha\in K_0(\ca)$, there exist $A,B\in\ca$ such that $\alpha=\widehat{A}-\widehat{B}\in K_0(\ca)$, and we set
	\begin{align}
		[K_\alpha]:=[K_A]* [K_B]^{-1},\qquad [K_\alpha^*]:=[K_A^*]*[K_B^*]^{-1}.
	\end{align}
	Then $[K_\alpha],[K_\alpha^*]$ are well defined in $\tMH(\ca)$; see \cite[\S3.4]{LP16}.

	\begin{theorem}[cf. \text{\cite[Theorem 3.20]{LP16}}]
		\label{theorem basis of modified hall algebra}
		We have the following.
		\begin{enumerate}
			\item
			$\tMH(\ca)$ has a basis given by
			$$\{[C_A\oplus C_B^*]* [K_\alpha]*[K_\beta^*]\mid [A],[B]\in\Iso(\ca),\alpha,\beta\in K_0(\ca)\}.$$
			\item
			For any $M=(\xymatrix{ M^0 \ar@<0.5ex>[r]^{d^0}& M^1 \ar@<0.5ex>[l]^{d^1}  })\in \cc_{\Z_2}(\ca)$, in  $\tMH(\ca)$ we have
			\begin{align*}
				[M]
				=&\sqq^{ -\langle\widehat{H^0(M)}-\widehat{H^1(M)},\;\widehat{\Im d^0}-\widehat{\Im d^1}\rangle}
				[C_{H^0(M)}^*\oplus C_{H^1(M)}]*[K_{\Im d^0}]* [K_{\Im d^1}^*].
			\end{align*}
		\end{enumerate}
	\end{theorem}


	
	
	\begin{lemma}
		[\text{\cite[Lemma 4.3]{LP16}}]
		\label{lem:KC}
		For any $\alpha,\beta\in K_0(\ca)$, $M\in \ca$, we have
		\begin{align}
			\big[[K_\alpha],[K_\beta]\big]=&\big[[K_\alpha],[K_\beta^*]\big]=\big[[K_\alpha^*],[K_\beta^*]\big]=0,
			\\
			[K_\alpha] *[C_M]=&\sqq^{(\alpha,\widehat{M})}[C_M]* [K_\alpha], \qquad \qquad [K_\alpha] *[C_M^*]= \sqq^{-(\alpha,\widehat{M})} [C_M^*]*[K_\alpha],
			\\
			[K_\alpha^*]*[C_M]=&\sqq^{-(\widehat{M},\alpha)} [C_M]*[K_\alpha^*],\qquad\quad [K_\alpha^*]*[C_M^*]=\sqq^{(\widehat{M},\alpha)}[C_M^*]*[K_\alpha^*].
		\end{align}
	\end{lemma}

	\section{Quantum Borcherds-Bozec algebras}
	\label{sec:Drinfeld}
	
	\subsection{Quantum Borcherds-Bozec algebras}
	\label{subsec: QBB}
	Let $\I$ be a finite or countably infinite index set. An (even) symmetric Borcherds-Cartan matrix is $A=(a_{ij})_{i,j\in \I}$ such that
	\begin{enumerate}
		\item $a_{ii}=2,0,-2,-4,\dots$;
		\item $a_{ij}=a_{ji}\in\Z_{\leq0}$ for $i\neq j$.
	\end{enumerate}
	Denote by $\fg$ the Borcherds-Bozec algebra associated to $A$.
	Let $\Pi=\{\alpha_i\mid i\in\I\}$ be the set of simple roots, and $\Z\I=\Z\Pi$ the root lattice.
	
	Let $\Ire:=\{i\in\I\mid a_{ii}=2\}$, $\Iim:=\{i\in\I\mid a_{ii}\leq0\}$. 
	Let $\I^\infty:=(\Ire\times\{1\})\cup (\Iim\times \Z_{>0})$. For simplicity, we shall often write $i$ instead of $(i,1)$ for $i\in\Ire$.
	
	For any $i\in\Ire$, the {\em simple reflection} $\bs_i:\Z^{\I}\rightarrow\Z^{\I}$ is defined to be $\bs_i(\alpha_j)=\alpha_j-a_{ij}\alpha_i$, for $j\in \I$.
	Denote the Weyl group by $W =\langle \bs_i\mid i\in \Ire\rangle$.
	
	Let $v$ be an indeterminant. Denote, for $r,m \in \N$,
	\[
	[r]=[r]_v=\frac{v^r-v^{-r}}{v-v^{-1}},
	\quad
	[r]!=[r]^!_v=\prod_{i=1}^r [i], \quad \qbinom{m}{r}=\qbinom{m}{r}_v =\frac{[m][m-1]\ldots [m-r+1]}{[r]!}.
	\]

	For any $i\in\I$, denote by
	\begin{align}
		v_{(i)}=v^{\frac{a_{ii}}{2}}.
	\end{align}
	
	\begin{definition}\label{definition of BC algebra}
		Let $\tau_{il}\in 1+v\Z_{\geq0}[[v]]$ for $ (i,l)\in\I^\infty$.
		The (Drinfeld double) quantum Borcherds-Bozec algebra $\tU:=\tU_v(\fg)$ associated to a Borcherds-Cartan matrix $A$ and the parameter $(\tau_{il})_{(i,l)\in\I^\infty}$ is the associative algebra over $\Q(v)$ generated by the elements $K_i^{\pm1}$, $(K_i')^{\pm1}$ ($i\in\I$), $e_{il}$, $f_{il}$ ($(i,l)\in\I^\infty$) with the following defining relations:
		\begin{align}
			\label{eq:QBB1}
			&K_iK_i^{-1}=K_i^{-1}K_i=1, \qquad K'_i(K'_i)^{-1}=(K'_i)^{-1}K'_i=1, 
			\\
			\label{eq:QBB2}
			&[K_i,K_j]=[K_i,K_j']=[K_i',K_j']=0, 
			\\
			\label{eq:QBB3}
			&K_ie_{jl}=v^{la_{ij}}e_{jl}K_i, \qquad
			K_if_{jl}=v^{-la_{ij}}f_{jl}K_i,
			\\
			\label{eq:QBB4}
			&K'_ie_{jl}=v^{-la_{ij}}e_{jl}K'_i, \qquad
			K'_if_{jl}=v^{la_{ij}}f_{jl}K'_i,
			\\
			\label{eq:QBB5}
			&e_{ik}f_{jl}-f_{jl}e_{ik}=0, \quad \text{ for }i\neq j,
			\\
			\label{eq:QBB6}
			&\sum_{\stackrel{m+r=k}{r+s=l}} v_{(i)}^{r(m-s)}\tau_{ir}e_{is}f_{im}(K_i')^r=\sum_{\stackrel{m+r=k}{r+s=l}} v_{(i)}^{-r(m-s)}\tau_{ir}f_{im}e_{is}(K_i)^r, \text{ for }(i,l),(i,k)\in\I^\infty,
			\\
			\label{eq:QBB7}
			&
			\sum_{k=0}^{1-la_{ij}}(-1)^k\qbinom{1-la_{ij}}{k} e_i^{1-la_{ij}-k}e_{jl}e_i^k=0,  \quad \text{ for }i\neq j,
			\\
			\label{eq:QBB8}
			&\sum_{k=0}^{1-la_{ij}}(-1)^k\qbinom{1-la_{ij}}{k} f_i^{1-la_{ij}-k}f_{jl}f_i^k=0,  \quad \text{ for }i\neq j.
		\end{align}
	\end{definition}

	For any $r\geq0$, define
	\begin{align}
		\varphi_r(t)=(1-t)(1-t^2)\cdots(1-t^r).
	\end{align}
	From now on, we shall assume that
	\begin{align}
		\label{eq:tauspec}
		\tau_{ir}=\frac{1}{\varphi_r(v^2)}, \qquad\text{ for any }(i,r)\in\I^\infty.
	\end{align}

The quantum Borcherds-Bozec algebra $\U:=\U_v(\fg)$  is defined to be the $\Q(v)$-algebra generated by the elements $K_i^{\pm1}$ ($i\in\I$), $e_{il}$, $f_{il}$ ($(i,l)\in\I^\infty$) subject to the defining relations \eqref{eq:QBB1}--\eqref{eq:QBB8} with $K_i'$ replaced by $K_i^{-1}$; see  \cite{B15,B16,FKKT20}.
	In particular, $\U$ is isomorphic to the quotient algebra of $\tU$ modulo the ideal $(K_iK_i'-1\mid i\in\I)$.

	Let $\tU^+$ be the subalgebra of $\tU$ generated by $e_{il}$ $((i,l)\in \I^\infty)$, $\tU^0$ be the subalgebra of $\tU$ generated by $K_i, K_i'$ $(i\in \I)$, and $\tU^-$ be the subalgebra of $\tU$ generated by $f_{il}$ $((i,l)\in \I^\infty)$, respectively.
	The subalgebras $\U^+$, $\U^0$ and $\U^-$ of $\U$ are defined similarly. Then both $\tU$ and $\U$ have triangular decompositions:
	\begin{align}
		\label{triandecomp2}
		\tU =\tU^+\otimes \tU^0\otimes\tU^-,
		\qquad
		\U &=\U^+\otimes \U^0\otimes\U^-.
	\end{align}
	Clearly, ${\U}^+\cong\tU^+$, $\U^-\cong \tU^-$, and $\U^0 \cong \tU^0/(K_i K_i' -1 \mid   i\in \I)$.
	

	Let $\omega$ be the involution of $\tU$ such that
	\begin{align*}
		\omega(K_i)=K_i',\qquad 	\omega(K'_i)=K_i,\qquad \omega(e_{il})=f_{il},\qquad \omega(f_{il})=e_{il},\quad \text{ for any }(i,l)\in\I^{\infty}.
	\end{align*}

Obviously, $\tU=\bigoplus_{\mu\in\Z\I} \tU_{\mu}$ is a graded algebra by setting $\deg K_i=0=\deg K_i'$, $\deg e_{il}=l\alpha_i$ and $\deg f_{il}=-l\alpha_i$.
In \cite{B15,B16}, Bozec introduced the {\em primitive generators} $\mathfrak{a}_{il}\in\tU^+_{l\alpha_i}$,  $\mathfrak{b}_{il}=\omega(a_{il})\in\tU^-_{-l\alpha_i}$; see e.g. \cite[Proposition 1.3]{FT21}. We do not need the precise definition of $\mathfrak{a}_{il}$ and  $\mathfrak{b}_{il}$ here, except the property
	\begin{align}
	\label{def:ail}
	\mathfrak{a}_{il}-e_{il}\in\Q(v)\langle e_{ik}\mid k<l\rangle, \text{ and }	\mathfrak{b}_{il}-f_{il}\in\Q(v)\langle f_{ik}\mid k<l\rangle.
\end{align}

	Given an $i\in\Ire$, as in \cite[Section 2]{FT21}, define the symmetry $T_{i}=T'_{i,1}:\tU\rightarrow \tU$ on the generators as follows:
	\begin{align*}
		T_i(\mathfrak{a}_i)=&v(K'_i)^{-1}\mathfrak{b}_i	,\qquad T_i(\mathfrak{b}_i)=v^{-1}\mathfrak{a}_iK_i^{-1},\\
		T_i(\mathfrak{a}_{jl})=&\sum_{r+s=-la_{ij}}(-1)^{r}v^{r}\mathfrak{a}_i^{(r)}\mathfrak{a}_{jl}\mathfrak{a}_i^{(s)}, \quad \text{ for }j\neq i,
		\\
		T_i(\mathfrak{b}_{jl})=&\sum_{r+s=-la_{ij}}(-1)^{r}v^{r}\mathfrak{b}_i^{(r)}{\mathfrak b}_{jl}\mathfrak{b}_i^{(s)},  \quad \text{ for }j\neq i,
		\\
		T_i(K_\mu)=&K_{\bs_i\mu},\qquad T_i(K'_\mu)=K'_{\bs_i\mu},\text{ for any }\mu\in\Z\I.
	\end{align*}
	Here we set
	$K_\mu:=\prod_{i\in\I} K_i^{l_i}$, and $K_\mu':=\prod_{i\in\I}(K_i')^{l_i}$, for any $\mu=\sum_{i\in\I} l_i\alpha_i\in\Z\I$.

	The symmetry $T''_{i,-1}:\tU\rightarrow \tU$ is defined such that
	$$T''_{i,-1}=\sigma T'_{i,1}\sigma,$$
	where  $\sigma:\tU\rightarrow \tU$ is the anti-involution such that
	\begin{align*}
		\sigma(e_{il})=e_{il},\qquad \sigma(f_{il})=f_{il},\qquad \sigma(K_i)=K_i',\quad \forall (i,l)\in\I^\infty.
	\end{align*}
	
	Indeed, in \cite{FT21}, they also define other two symmetries $T'_{i,-1}$ and $T''_{i,1}$. These symmetries are algebra isomorphisms and satisfy the braid relations associated to a Borcherds-Cartan datum.
	
	In fact, the quantum Borcherds-Bozec algebra $\U$ in \cite{FKKT20} is slightly different to the one  defined in \cite{FT21},  which is denoted by $\U'$ to avoid confusions. We also use $\tU'$ to denote the Drinfeld double version of $\U'$.  However, they coincide up to the algebra isomorphism $\phi: \tU\rightarrow \tU'$ over $\Q$, which sends
 $$v\mapsto v^{-1},\quad e_{il}\mapsto e_{il},\quad f_{il}\mapsto f_{il},\quad K_i\mapsto K_i', \quad K_i'\mapsto K_i.$$
	Conjugating with $\phi$, the braid group actions $T'_{i,e}$ and $T''_{i,e}$  coincide with $T'_{i,-e}$ and $T''_{i,-e}$ defined in \cite{FT21}, respectively .

	\subsection{Semi-derived Ringel-Hall algebras of quivers}
	Let $Q=(Q_0=\I,\Omega)$ be a quiver, where $\I$ is the set of vertices and $\Omega$ is  the set of arrows. Define $g_i$ the number of loops at $i$ and $n_{ij}$ the number of arrows from $i$ to $j$ if $i\neq j$. Let $A=(a_{ij})_{i,j\in \I}$ be the symmetric Borcherds-Cartan matrix associated to $Q$, i.e.,
	\begin{align}
		\label{eq:BC-Q}
		a_{ij}=\begin{cases}2-2g_i& \text{ if }i=j, \\
		-n_{ij}-n_{ji} & \text{ if }i\neq j.\end{cases}
		\end{align}
	
	A representation of $Q$ over $\bfk$ is $(M_i,M_\alpha)_{i\in\I,\alpha\in\Omega}$ by associating a finite-dimensional $\bfk$-linear space $M_i$ to each vertex $i\in\I$, and a  $\bfk$-linear map $M_\alpha: M_i\rightarrow M_j$ to each arrow $\alpha:i\rightarrow j\in\Omega$.
	A representation $(M_i,M_\alpha)$ of $Q$ is called nilpotent if $M_{\alpha_r}\cdots M_{\alpha_{2}}M_{\alpha_1}$ is nilpotent for any cyclic paths $\alpha_r\cdots \alpha_2\alpha_1$.
	
	Let $\rep_\bfk^{\rm nil}(Q)$ be the category of finite-dimensional nilpotent representations of $Q$ over $\bfk$. Then $\rep_\bfk^{\rm nil}(Q)$ is a hereditary abelian category.
	Let $\langle\cdot,\cdot\rangle_Q$ be the Euler form of $Q$. Define
	$$(x,y)_Q=\langle x,y\rangle_Q+\langle y,x\rangle_Q.$$
	Let $S_i$ be the nilpotent simple module supported at $i\in\I$.
	Then $(S_i,S_j)_Q=a_{ij}$ for any $i,j\in\I$; see e.g. \cite[Chapter 1, Example 1.26]{Ki16}.
	
	Denote by $\sqq_{(i)}=\sqq^{\frac{a_{ii}}{2}}=\sqq^{1-g_i}$ for $i\in\I$.
	
 The quantum Borcherds-Bozec algebra $\tU$ associated to the quiver $Q$ is the one constructed by using the Borcherds-Cartan matrix $A$ associated to $Q$ and the parameter $(\tau_{il})_{(i,l)\in\I^\infty}$ in \eqref{eq:tauspec}.

	\begin{theorem}[\cite{Lu22}]
		\label{thm:iso}
		Let $Q$ be an arbitrary quiver, and $\tU$ be its quantum Borcherds-Bozec algebra. Then there is an injective morphism of algebras $\Psi: \tU\rightarrow \tMH(\rep^{\rm nil}_\bfk(Q))$ such that
		\begin{align*}
			\Psi(K_i)=&[K_{S_i}],\quad \Psi(K_i')=[K_{S_i}^*],\\
			\Psi(e_{il})=&(-1)^l\sqq^{l^2-l}[\![C_{S_i^{\oplus l}}]\!],\quad \Psi(f_{il})=(-1)^l\sqq^{l^2-l}[\![C^*_{S_i^{\oplus l}}]\!]
		\end{align*}
		for any $i\in\I$, $(i,l)\in\I^\infty$.
	\end{theorem}
	
	\begin{proof}
		The morphism given here is slightly different to \cite[Theorem 3.4]{Lu22}. However, by  \cite[Theorem 3.4]{Lu22}, one only needs to check $\Psi$ preserves \eqref{eq:QBB6}, which follows directly by using \cite[Lemma  3.3]{Lu22}.	
	\end{proof}

	\subsection{Reflection functors}
	
	For any vertex $\ell$ of $Q$, let $\bs_\ell Q$ be the quiver constructed from $Q$ by reversing all arrows starting or ending to $\ell$.
	
	Fix a sink vertex $\ell$ in $Q$. Let
	\begin{align}
		F_\ell^+:\rep_\bfk^{\rm nil}(Q)\longrightarrow \rep_\bfk^{\rm nil}(\bs_\ell Q)
	\end{align}
	be the BGP reflection functor; see \cite{BGP73}. Then $F_\ell^+$ induces a functor  (denoted by the same notation):
	\begin{align}
		F_\ell^+:\cc_{\Z_2}(\rep_\bfk^{\rm nil}(Q))\longrightarrow \cc_{\Z_2}(\rep_\bfk^{\rm nil}(\bs_\ell Q)).
	\end{align}

	Let $\ct=\ct_\ell:=\{X\in\rep_\bfk^{\rm nil}(Q)\mid \Hom_{\bfk Q}(X,S_\ell)=0\}$, and let $\cf=\cf_\ell$ be the extension closed subcategory of $\rep_\bfk^{\rm nil}( Q)$ generated by $S_\ell$. Then $(\ct,\cf)$ is a torsion pair.
	Then $(\cc_{\Z_2}(\ct),\cc_{\Z_2}(\cf))$ is also a torsion pair of $\cc_{\Z_2}(\rep_\bfk^{\rm nil}(Q))$; see \cite[Lemma 2.7]{LW22}. In particular, for any $M\in \cc_{\Z_2}(\rep_\bfk^{\rm nil}(Q))$, there exists a short exact sequence
	\begin{align}
		\label{resol}
		0 \longrightarrow M \longrightarrow X_M \longrightarrow T_M \longrightarrow 0
	\end{align}
	with $X_M,T_M\in \cc_{\Z_2}(\ct)$ and $T_M\in\cc_{\Z_2,ac}(\rep_\bfk^{\rm nil}(Q))$; see \cite[Lemma 2.7]{LW22}~(b);  cf. \cite[Proposition 5.8]{LP16}.

	\begin{lemma}[\text{cf. \cite[Theorem 2.11]{LW22}}]
		\label{lem:Upsilon}
		Let $\ell$ be a sink of $Q$. Then we have an isomorphism of algebras:
		\begin{align}
			\label{eqn:reflection functor 1}
			\Gamma_{\ell}:\cs\cd\widetilde{\ch}(\bfk Q)&\stackrel{\cong}{\longrightarrow} \cs\cd\widetilde{\ch}(\bfk \bs_\ell Q)\\\notag
			[M]&\mapsto \sqq^{\langle  \res T_M,\res M\rangle_Q} q^{-\langle T_M,M\rangle}  [F_\ell^+(T_M)]^{-1}\diamond [F_\ell^+(X_M)],
		\end{align}
		where $M\in\cc_{\Z_2}(\rep_\bfk^{\rm nil}(Q))$ and $X_M,T_M\in \cc_{\Z_2}(\ct)$ ($T_M$ acyclic) fit \eqref{resol}.
	\end{lemma}

\begin{proof}
		The $\imath$quiver algebra and its module category are considered in \cite{LW19a}, and the category $\cc_{\Z_2}(\rep_\bfk^{\rm nil}(Q))$ can be viewed to be the module category of an $\imath$quiver algebra of diagonal type; see \cite[Example 2.10]{LW19a}. So the statement follows from \cite[Theorem 2.11]{LW22}.
\end{proof}

	Similarly, for any source $\ell$ of $Q$, we have an isomorphism of algebras
	$$\Gamma_\ell^-:\cs\cd\widetilde{\ch}(\bfk Q)\stackrel{\cong}{\longrightarrow} \cs\cd\widetilde{\ch}(\bfk \bs_\ell Q).$$ 
	In particular, for a quiver $Q$ with $\ell$ a sink, we have
	$\Gamma_{\ell}:\cs\cd\widetilde{\ch}(\bfk Q){\rightarrow} \cs\cd\widetilde{\ch}(\bfk \bs_\ell Q)$ and $\Gamma_\ell^-:\cs\cd\widetilde{\ch}(\bfk\bs_\ell Q){\longrightarrow} \cs\cd\widetilde{\ch}(\bfk Q)$ are mutually inverses to each other.
	
	\begin{proposition}[\text{cf. \cite[Proposition 2.12]{LW22}}]
		\label{prop:reflection}
		Let $Q$ be a quiver, and $\ell\in Q_0$ be a sink.
		Then the isomorphism $\Gamma_{\ell}:\cs\cd\widetilde{\ch}(\bfk Q)\stackrel{\cong}{\longrightarrow} \cs\cd\widetilde{\ch}(\bfk \bs_\ell Q)$ sends
		\begin{align}
			\Gamma_{\ell}([M])&= [F_{\ell}^+(M)], \quad \forall M\in\cc_{\Z_2}(\ct),
			\label{eqn:reflection 1}
			\\
			\Gamma_{\ell}([C_{S_\ell}])&=
			\sqq[K^*_{S_\ell}]^{-1}* [C^*_{S_{ \ell}}],\qquad \Gamma_{\ell}([C^*_{S_\ell}])=
			\sqq[K_{S_\ell}]^{-1}* [C_{S_{ \ell}}],
			\label{eqn:reflection 2}
			\\
			\Gamma_{\ell}([K_\alpha])&= [K_{\bs_{\ell}\alpha}], \qquad \Gamma_{\ell}([K^*_\alpha])= [K^*_{\bs_{\ell}\alpha}],
			\quad \forall\alpha\in K_0(\rep_\bfk^{\rm nil}(Q)).\label{eqn:reflection 4}
		\end{align}
	\end{proposition}

	\begin{proof}
By the proof of Lemma \ref{lem:Upsilon}, the statement follows from \cite[Proposition 2.12]{LW22}.
	\end{proof}

	\begin{theorem}
		\label{thm:main1}
		Let $Q$ be a quiver, and $\ell\in \Ire$ be a sink. Then we have the following commutative diagram:
		\begin{equation}
			\label{eq:commdiag}
			\xymatrix{\tU|_{v=\sqq}\ar[r]^{T_\ell} \ar[d]^{\Psi_Q}&  \tU|_{v=\sqq}\ar[d]^{\Psi_{\bs_\ell Q}} \\
				\cs\cd\widetilde{\ch}(\bfk Q)\ar[r]^{\Gamma_\ell} & \cs\cd\widetilde{\ch}(\bfk \bs_\ell Q)	}
		\end{equation}
	\end{theorem}
	
	The proof of Theorem
	\ref{thm:main1} shall be given in the following subsections.

	Similarly, we have have the following result.
	
	\begin{proposition}
		\label{prop:source}
		Let $Q$ be a quiver. For any source $\ell\in \Ire$, one can obtain the following commutative diagram:
		\begin{equation}
			\label{eq:commdiag-new}
			\xymatrix{\tU|_{v=\sqq}\ar[r]^{T''_{\ell,-1}} \ar[d]^{\Psi_Q}&  \tU|_{v=\sqq}\ar[d]^{\Psi_{\bs_\ell Q}} \\
				\cs\cd\widetilde{\ch}(\bfk Q)\ar[r]^{\Gamma^-_\ell} & \cs\cd\widetilde{\ch}(\bfk \bs_\ell Q)	}
		\end{equation}
	\end{proposition}
	
	\begin{corollary}
		For any $\ell\in \Ire$, $T_{\ell}=T'_{\ell,1}$ and $T''_{\ell,-1}$
		are algebra endomorphisms of $\tU$,  which are mutually inverses to each other.
	\end{corollary}
	
	\begin{proof}
		By the commutative diagram obtained in	Theorem \ref{thm:main1} and Proposition \ref{prop:source}, one can see that $T'_{\ell,1},T''_{\ell,-1}:
		\tU\rightarrow\tU$ are algebra homomorphisms.
		
		For any sink $\ell\in\Ire$, we have
		$\Gamma_\ell: \cs\cd\widetilde{\ch}(\bfk Q)\rightarrow \cs\cd\widetilde{\ch}(\bfk \bs_\ell Q)$ and $\Gamma^-_\ell:\cs\cd\widetilde{\ch}(\bfk\bs_\ell Q)\rightarrow \cs\cd\widetilde{\ch}(\bfk Q)$ are mutually inverses to each other.
		Then the desired result follows from the commutative diagrams \eqref{eq:commdiag} and \eqref{eq:commdiag-new}.
	\end{proof}

	\begin{remark}
		The actions of $T'_{i,1}, T''_{i,-1}$ on $\tU$ factor through the quotient $\U =\tU \big/ ( K_i' K_i-1\mid i\in\I )$ to the corresponding automorphisms on $\U$, and the formulas \eqref{eqn:reflection 1}--\eqref{eqn:reflection 4} are then reduced to  formulas in \cite{FT21} upon the identification $K_i' =K_i^{-1}$.
	\end{remark}
	
	\subsection{A building block}
	Let $Q$ be a quiver, and $\ell\in Q_0$ be a sink. For any vertex $j\neq \ell$ and any $\bfk Q$-module $N$ supported on $j$, we can (shall always) view $N$ to be a $\bfk\bs_\ell Q$-module naturally.
	
	\begin{proposition}
		\label{prop:building}
		Let $Q$ be a quiver, and $\ell\in Q_0$ be a sink. For any vertex $j\neq \ell$, let $N$ be a nilpotent $\bfk Q$-module of dimension $l$ supported on $j$. Then in $\cs\cd\widetilde{\ch}(\bfk \bs_\ell Q)$ we have
		\begin{align}
			\label{eq:RefHallE}
			\Gamma_\ell([C_N])=(1-q)^{la_{\ell,j}}\sum_{r+s=-la_{\ell,j}}(-1)^r \sqq^{r}[C_{S_\ell}]^{(r)}*[C_N]*[C_{S_\ell}]^{(s)},
			\\
			\label{eq:RefHallF}
			\Gamma_\ell([C^*_N])=(1-q)^{la_{\ell,j}}\sum_{r+s=-la_{\ell,j}}(-1)^r \sqq^{r}[C^*_{S_\ell}]^{(r)}*[C^*_N]*[C^*_{S_\ell}]^{(s)}.
		\end{align}
	\end{proposition}

	\begin{proof}
		We only prove \eqref{eq:RefHallE} since
		\eqref{eq:RefHallF} is similar.
		
		Denote by $a=-a_{\ell,j}$. For any $r,s\geq0$ with $r+s=la$, we have
		$$\langle C_{S_\ell^{\oplus r}}, C_N\rangle=-rla,\quad \langle C_{S_\ell^{\oplus r}}, C_{S_\ell^{\oplus s}}\rangle=rs.$$
		Hence,
		\begin{align}\label{quantum relation for module} &[C_{S_\ell^{\oplus r}}]*[C_N]*[C_{S_\ell^{\oplus s}}]\\\notag
			=&[C_{S_\ell^{\oplus r}}]*[C_{N\oplus S_\ell^{\oplus s}}]\\\notag
			=&\sqq^{rs-rla}\sum\limits_{[C_M]\in I_{r,s,l}}\frac{|\Ext^1(C_{S_\ell^{\oplus r}}, C_{N\oplus S_\ell^{\oplus s}})_{C_M}|}{|\Hom(C_{S_\ell^{\oplus r}}, C_{N\oplus S_\ell^{\oplus s}})|}[C_M]\\\notag
			=&\sqq^{rs-rla}\sum\limits_{[C_M]\in I_{r,s,l}}F_{C_{S_\ell^{\oplus r}}, C_{N\oplus S_\ell^{\oplus s}}}^{C_M}\cdot \frac{|\Aut(C_{S_\ell^{\oplus r}})|\cdot |\Aut(C_{N\oplus S_\ell^{\oplus s}})|}{|\Aut(C_M)|}[C_M],
		\end{align}
		where $$I_{r,s,l}=\{[C_M]\;|\;\exists L\subseteq M, \;s.t.\; L\cong N\oplus S_\ell^{\oplus s},\; M/L\cong S_\ell^{\oplus r}\}.$$

	Assume $M=(M_i,M_\alpha)$.	For any $[C_M]\in I_{r,s,l}$, we have 
$\dim M_j=lae_{\ell}+le_{j}$.  
	Denote by $\alpha_i:\ell\rightarrow j$ ($1\leq i\leq a$) the arrows between $\ell$ and  $j$ in $\bs_\ell Q$.
		Let $U_{M}=\bigcap\limits_{i=1}^{a}\ker M_{\alpha_i}$ and $u_{M}=\dim U_{M}$. Then $u_{M}\geq s$, and
		$$F_{C_{S_\ell^{\oplus r}}, C_{N\oplus S_\ell^{\oplus s}}}^{C_M}=|{\rm{Gr}}(s, u_M)|=\sqq^{(u_M-s)s}\qbinom{u_M}{s}.$$
		
		Recall from \cite[(5.2)]{LW20} that for any $r>0$, $$[C_{S_\ell}]^{(r)}=\sqq^{-\frac{r(r-1)}{2}}\cdot
		\frac{[C_{S_\ell^{\oplus r}}]}{[r]^!_\sqq}.$$
		Moreover,
		\begin{align}  \sqq^{-\frac{r(r-1)}{2}}\cdot\frac{|\Aut(C_{S_\ell^{\oplus r}})|}{[r]^!_\sqq}
			=&\sqq^{-\frac{r(r-1)}{2}}\cdot\frac{(q^r-1)(q^r-q)\cdots(q^r-q^{r-1})}{[r]_\sqq[r-1]_\sqq\cdots [1]_\sqq}\\\notag
			=&\sqq^{-\frac{r(r-1)}{2}}(\sqq-\sqq^{-1})^r\cdot\prod_{i=1}^{r-1} q^{i}\cdot \prod_{i=1}^{r}\sqq^{i} \\\notag
			=&(q-1)^r\sqq^{r(r-1)}.
		\end{align}
		Then it follows from \eqref{quantum relation for module} that
		\begin{align}\label{quantum relation for divided power}
			&[C_{S_\ell}]^{(r)}*[C_N]*[C_{S_\ell}]^{(s)}\\\notag
			=&\sqq^{rs-rla}\sum\limits_{[C_M]\in I_{r,s,l}}
			(q-1)^r\sqq^{r(r-1)}\cdot(q-1)^s\sqq^{s(s-1)}\cdot\sqq^{(u_M-s)s}\qbinom{u_M}{s}_\sqq\cdot \frac{|\Aut(C_{N})|}{|\Aut(C_M)|}[C_M]\\\notag
			=&\sum\limits_{[C_M]\in I_{r,s,l}}
			(q-1)^{la}\sqq^{rs-rla+r(r-1)+s(s-1)+(u_M-s)s}\qbinom{u_M}{s}_\sqq\cdot\frac{|\Aut(C_{N})|}{|\Aut(C_M)|}[C_M].
		\end{align}
		To sum up, we obtain that
		\begin{align}
			\label{eq:CCC}
			&\sum_{r+s=la}(-1)^r \sqq^{r}[C_{S_\ell}]^{(r)}*[C_N]*[C_{S_\ell}]^{(s)}\\\notag
			=&\sum\limits_{[C_M]\in I_{r,s,l}}
			\sum_{r+s=la}(-1)^r \sqq^{r}\cdot
			(q-1)^{la}\sqq^{rs-rla+r(r-1)+s(s-1)+(u_M-s)s}\qbinom{u_M}{s}_\sqq\frac{|\Aut(C_{N})|}{|\Aut(C_M)|}
\cdot[C_M]\\\notag
=&(1-q)^{la}\sum\limits_{[C_M]\in I_{r,s,l}}\frac{|\Aut(C_{N})|}{|\Aut(C_M)|}
\cdot\big(\sum_{s\geq0}(-1)^s\sqq^{(u_M-1)s}\qbinom{u_M}{s}_\sqq\big)\cdot[C_M].
		\end{align}
	The coefficient of $[C_M]$ in \eqref{eq:CCC}  is zero if $u_M>0$ by using the standard $v$-binomial identity \cite[1.3.1 (c)]{Lus93}.
	
Note that $N\in\ct$.	We claim $M\cong F^+_\ell(N)$ if $u_M=0$.
		If the claim holds, then \eqref{eq:CCC} equals to $(1-q)^{la}\Gamma_\ell([C_N])$   since the  restriction of $F^+_\ell$ to $\ct$ is an embedding.
		
	It remains to prove the claim. Let $N=(N_i,N_\alpha)$, and then only $N_j\neq 0$.
	Let $F^+_\ell(N)=(Y_i,Y_\alpha)$.  Then  $Y_j=N_j$, $Y_\ell=N_j^{\oplus a}$, and $Y_{\alpha}=N_\alpha$ for any loop $\alpha$ at $j$, $Y_{\alpha_k}:N_j^{\oplus a}\rightarrow N_j$
is the $k$-th projection for $1\leq k\leq a$.

For $M$ with $u_M=0$, we have $M_j=N_j$ and $\dim M_\ell=la=\dim Y_\ell$.
Let $f=(f_i)_{i\in \I}:M\rightarrow F^+_\ell(N)$ be the morphism defined by
\begin{align*}
	f_i=\begin{cases} {\rm Id} & \text{ if }i=j,
	\\
	(\iota_k M_{\alpha_k})_k:M_\ell \rightarrow N_j^{\oplus a} & \text{ if }i=\ell,
	\\
	0 &\text{ otherwise},
  \end{cases}
\end{align*}
where $\iota_k: N_j\rightarrow N_j^{\oplus a}$ is the $k$-th injection. Obviously, $f$ is well defined, and $\ker f_\ell=U_M=0$. So $f$ is injective, which implies that $f$ is an isomorphism by considering the dimensions. The claim is proved.		
	\end{proof}
	
		\begin{remark}
	In case that $Q$ admits no loops, i.e., in the setting of quantum groups, the observation in Proposition \ref{prop:building} can be used to deduce the action of braid group symmetries on dividied powers (see \cite[Page 288]{Lus93}).
	\end{remark}
	
	\subsection{Proof of Theorem \ref{thm:main1}}
	In this subsection, we complete the proof of Theorem \ref{thm:main1}.
	
	\begin{proof}[Proof of Theorem \ref{thm:main1}]
		Note $\ell\in\I^{\rm re}$. So $\mathfrak{a}_\ell=e_\ell$ and $\mathfrak{b}_\ell=f_\ell$.
		
		By definition, we have
		
		\begin{align*}
			\Psi_{\bs_\ell Q} T_\ell(K_\mu)=\Gamma_\ell\Psi_Q(K_\mu),\qquad \Psi_{\bs_\ell Q} T_\ell(K'_\mu)=\Gamma_\ell\Psi_Q(K'_\mu),\quad \forall \mu\in\Z\I.
		\end{align*}
		Furthermore, we have
		\begin{align*}
			\Gamma_\ell\Psi_Q(\mathfrak{a}_\ell)&=\Gamma_\ell\Psi_Q(e_\ell)
			=\Gamma_\ell([\![C_{S_\ell}]\!])
			=\sqq[K^*_{S_\ell}]^{-1}* [\![C^*_{S_{ \ell}}]\!];
		\end{align*}
		hence
		\begin{align*}
			\Psi_{\bs_\ell Q} T_\ell(\mathfrak{a}_\ell)&=\sqq\Psi_{\bs_\ell Q} ((K_\ell')^{-1}\mathfrak{b}_\ell)
			=\sqq[K^*_{S_\ell}]^{-1}* [C^*_{S_{ \ell}}]=\Gamma_\ell\Psi_Q(\mathfrak{a}_\ell).
		\end{align*}
		
		Similarly,
		\begin{align*}
			\Gamma_\ell\Psi_Q(\mathfrak{b}_\ell)&=\Gamma_\ell\Psi_Q(f_\ell)
			=\Gamma_\ell([\![C^*_{S_\ell}]\!])
			=\sqq[K_{S_\ell}]^{-1}* [\![C_{S_{ \ell}}]\!]=\sqq^{-1} [\![C_{S_{ \ell}}]\!]*[K_{S_\ell}]^{-1};
		\end{align*}
		then
		\begin{align*}
			\Psi_{\bs_\ell Q} T_\ell(\mathfrak{b}_\ell)&=\sqq^{-1}\Psi_{\bs_\ell Q} (\mathfrak{a}_\ell K_\ell^{-1})
			= \sqq^{-1}[C_{S_{ \ell}}]*[K_{S_\ell}]^{-1}=\Gamma_\ell\Psi_Q(\mathfrak{b}_\ell).
		\end{align*}

		Using \eqref{def:ail} and the definition of $\Psi_Q$ in Theorem \ref{thm:iso}, one can see
		$$\Psi_Q(\mathfrak{a}_{jl})=\sum_{k\in\mathcal{K}} f_k(\sqq) [M_k]$$
		summed over a finite set $\mathcal{K}$, with  $f_k(\sqq)\in\Q(\sqq)$ and $M_k$ supported on $j$ with  dimension $l$ for any $k$.
		
		By Proposition \ref{prop:building}, we have
		\begin{align*}
			\Gamma_\ell\Psi_Q(\mathfrak{a}_{jl}) =&(1-q)^{la_{\ell,j}}\sum_{r+s=-la_{\ell,j}}(-1)^r \sqq^{r}[C_{S_\ell}]^{(r)}*\Psi_Q(\mathfrak{a}_{jl})*[C_{S_\ell}]^{(s)};
		\end{align*}
		and then
		\begin{align*}
			\Psi_{\bs_\ell Q} T_\ell(\mathfrak{a}_{jl})&=\Psi_{\bs_\ell Q} \big(\sum_{r+s=-la_{\ell,j}}(-1)^{r}v^{r}\mathfrak{a}_\ell^{(r)}\mathfrak{a}_{jl}\mathfrak{a}_\ell^{(s)}\big)
			\\
			&= (-1)^{la_{\ell,j}}\sum_{r+s=-la_{\ell,j}}(-1)^r \sqq^{r}[\![C_{S_\ell}]\!]^{(r)}*\Psi_{\bs_\ell Q}(\mathfrak{a}_{jl})*[\![C_{S_\ell}]\!]^{(s)}
			\\
			&= \Gamma_\ell\Psi_Q(\mathfrak{a}_{jl}).
		\end{align*}
		Similarly, one can check that
		\begin{align*}
			\Psi_{\bs_\ell Q} T_\ell(\mathfrak{b}_{jl})
			&= \Gamma_\ell\Psi_Q(\mathfrak{b}_{jl}).
		\end{align*}
		
		Therefore, the commutative diagram \eqref{eq:commdiag} follows.
	\end{proof}

	\begin{corollary}
		For any $i\in \I^{\rm re}$ and $(j,l)\in\I^{\infty}$ with $j\neq i$, we have
		\begin{align}
			T_i(e_{jl})=&\sum_{r+s=-la_{ij}}(-1)^{r}v^{r}e_i^{(r)}e_{jl}e_i^{(s)},
			\\
			T_i(f_{jl})=&\sum_{r+s=-la_{ij}}(-1)^{r}v^{r}f_i^{(r)}f_{jl}f_i^{(s)}.
		\end{align}	
	\end{corollary}

	\begin{proof}
		It follows from Theorem \ref{thm:main1} and Proposition \ref{prop:building}.	
	\end{proof}	
	
	\section{Quantum generalized Kac-Moody algebras}
	\label{sec:Quantum generalized}

	In this section, we use semi-derived Ringel-Hall algebras of quivers with loops to establish the braid group actions of the quantum generalized Kac-Moody algebras briefly.
	
	Let $A=(a_{ij})_{i,j\in\I}$ be a symmetric Borcherds-Cartan matrix as defined in \S\ref{subsec: QBB}. Assume that we are given a collection of positive integers $\bm=(m_i)_{i\in\I}$ with $m_i=1$ whenever $i\in\Ire$, called the \emph{charge} of $A$. Let $\fg_{A,\bm}$ be the generalized Kac-Moody algebra (also called Borcherds algebra).
	
	The (Drinfeld double) quantum generalized Kac-Moody algebra $\tUB:=\tUB_v(\fg_{A,\bm})$ is the $\Q(v)$-algebra generated by the elements $K_i^{\pm1}$, $(K_i')^{\pm1}$, $E_{ik}$, $F_{ik}$ for $i\in\I$, $k=1,\dots,m_i$ subject to \eqref{eq:QBB1}--\eqref{eq:QBB2} and
	\begin{align}
		\label{eq:QB3}
		&K_iE_{jl}=v^{a_{ij}}E_{jl}K_i, \qquad
		K_iF_{jl}=v^{-a_{ij}}F_{jl}K_i,
		\\
		\label{eq:QB4}
		&K'_iE_{jl}=v^{-a_{ij}}E_{jl}K'_i, \qquad
		K'_iF_{jl}=v^{a_{ij}}F_{jl}K'_i,
		\\
		\label{eq:QB5}
		&E_{ik}F_{jl}-F_{jl}E_{ik}=\delta_{ij}\delta_{kl}\frac{K_i-K_i'}{v-v^{-1}},
		\\
		\label{eq:QB7}
		&
		\sum_{n=0}^{1-a_{ij}}(-1)^n\qbinom{1-a_{ij}}{n} E_{i1}^{1-a_{ij}-n}E_{jl}E_{i1}^n=0, 	\quad \text{ for } i\neq j,
		\\
		\label{eq:QB8}
		&\sum_{n=0}^{1-a_{ij}}(-1)^n\qbinom{1-a_{ij}}{n} F_{i1}^{1-a_{ij}-n}F_{jl}F_{i1}^n=0,
		\quad \text{ for } i\neq j.
	\end{align}
	
	Similar to quantum groups, there exists a $\Q(v)$-algebra anti-involution $\sigma:\tUB\rightarrow \tUB$ such that
	\begin{align*}
		\sigma(E_{il})=E_{il},\quad \sigma(F_{il})=F_{il},\quad \sigma(K_i)= K'_{i},
		\quad \forall (i,l)\in \I^\infty.
	\end{align*}
	There exists a $\Q$-algebra automorphism $ \psi: \tUB\rightarrow \tUB$ (called \emph{bar involution}) such that
	\[
	\psi(v)=v^{-1}, \quad
	\psi(K_i)=K'_{i}, \quad
	\psi (E_{il})=E_{il},  \quad \psi(F_{il})=F_{il},\quad
	\forall (i,l)\in \I^\infty.
	\]

	Note that $K_iK_i'$ are central in $\tUB$ for all $i\in\I$. The quantum generalized Kac-Moody algebra $\UB:=\UB_v(\fg_{A,\bm})$ is generated by the elements $K_i^{\pm1}$, $E_{ik}$, $F_{ik}$ for $i\in\I$, $k=1,\dots,m_i$ subject to \eqref{eq:QBB1}--\eqref{eq:QBB2}, \eqref{eq:QB3}--\eqref{eq:QB8} with $K_i'$ replaced by $K_i^{-1}$.
	In particular, the quantum generalized Kac-Moody algebra $\UB$ is isomorphic to the quotient algebra of $\tUB$ modulo the ideal $(K_iK_i'-1\mid i\in\I)$.
	
	Let $\tUB^+$ be the subalgebra of $\tUB$ generated by $E_{ik}$ $(i\in \I,k=1,\dots,m_i)$, $\tUB^0$ be the subalgebra of $\tUB$ generated by $K_i, K_i'$ $(i\in \I)$, and $\tUB^-$ be the subalgebra of $\tUB$ generated by $F_{ik}$ $(i\in \I,k=1,\dots,m_i)$, respectively.
	The subalgebras $\UB^+$, $\UB^0$ and $\UB^-$ of $\UB$ are defined similarly. Then both $\tUB$ and $\UB$ have triangular decompositions:
	\begin{align}
		\label{triandecomp}
		\tUB =\tUB^+\otimes \tUB^0\otimes\tUB^-,
		\qquad
		\UB &=\UB^+\otimes \UB^0\otimes\UB^-.
	\end{align}
	Clearly, ${\UB}^+\cong\tUB^+$, $\UB^-\cong \tUB^-$, and $\UB^0 \cong \tUB^0/(K_i K_i' -1 \mid   i\in \I)$.

	Let $Q=(\I,\Omega)$ be a quiver. Let $\rep_\bfk(Q)$ be the category of finite-dimensional representations of $Q$ over $\bfk$. Then $\rep_\bfk(Q)$ is a hereditary abelian category. We assume that the charge $\bm=(m_i)_{i\in\I}$ satisfies $m_i\leq |\bfk^{g_i}|=q^{g_i}$ for all $i\in\I$.
	If $i\in\Ire$ then there exists a unique simple representstion $S_i\in\rep_\bfk(Q)$ supported at $i$. If $i\in\Iim$ then the set of simple representations supported at $i$ is in bijection with $\bfk^{g_i}$: if $\sigma_1,\dots, \sigma_{g_i}$ denote the simple loops at $i$ then to $\underline{\lambda}=(\lambda_1,\dots,\lambda_{g_i})$ corresponds the simple representation $S_i(\underline{\lambda})=(W_j,x_\alpha)$ with $\dim_\bfk W_j=\delta_{ij}$ and $x_{\alpha_l}=\lambda_l$ for $l=1,\dots,g_i$.

	Following \cite{KS06}, we choose $\underline{\lambda}_i^{(l)}\in\bfk^{g_i}$ for $l=1,\dots,m_i$ in such a way that
	$\underline{\lambda}_i^{(l)}\neq \underline{\lambda}_i^{(l')}$ for $l\neq l'$. We set $S_{il}:=S_{i}(\underline{\lambda}_i^{(l)})$ for $i\in\Iim$ and $l=1,\dots, m_i$; and set $S_{i1}:=S_{i}$ for $i\in\Ire$.

	In order to realize the quantum generalized Kac-Moody algebras $\tUB$, we modify the definition of semi-derived Ringel-Hall algebras  slightly. The algebra $\cs\widetilde{\ch}(\rep_\bfk (Q))$ is constructed in the same way as $\tMH(\rep_\bfk (Q))$ but with the ideal $I_{\Z_2}$ of $\ch(\cc_{\Z_2}(\rep_\bfk(Q)))$ generated by
	\begin{align}
		&\{[L]-[K\oplus M]\mid \exists \text{ a short exact sequence }0 \longrightarrow K \longrightarrow L \longrightarrow M \longrightarrow 0 \}
		\\
		\notag
		\bigcup&\{ [K_M]-\prod_{i\in\I}[K_{S_i}]^{a_i},\quad [K^*_M]-\prod_{i\in\I}[K^*_{S_i}]^{a_i}\mid M\in\rep_\bfk(Q) \text{ with }\dimv M=(a_i)_{i\in\I} \}.
	\end{align}
	In fact, $\cs\widetilde{\ch}(\rep_\bfk (Q))$ is isomorphic to the  quotient algebra of $\tMH(\rep_\bfk (Q))$ modulo the ideal generated by
	\begin{align}
		\label{eq:ideal2}
		[K_M]-\prod_{i\in\I}[K_{S_i}]^{a_i},\quad [K^*_M]-\prod_{i\in\I}[K^*_{S_i}]^{a_i},\text{ for any }M\in\rep_\bfk(Q) \text{ with }\dimv M=(a_i)_{i\in\I}.
	\end{align}

	\begin{theorem}[\cite{Lu22}] Let $Q$ be an arbitrary quiver.
		There exists an injective morphism of algebras $\Psi=\Psi_Q: \tUB|_{v=\sqq}\rightarrow \cs\widetilde{\ch}(\rep_\bfk (Q))$ such that
		\begin{align}
			\Psi(K_i)=[K_{S_i}],\quad \Psi(K_i')=[K_{S_i}^*],\quad \Psi(E_{il})=\frac{1}{q-1}[C_{S_{il}}],\quad \Psi(F_{il})=\frac{-\sqq}{q-1}[C^*_{S_{il}}]
		\end{align}
		for any $i\in\I$.
	\end{theorem}
	
	Given an $i\in\Ire$, as in \cite[Chapter 37]{Lus93} (see also \cite[Section 2]{FT21}), define the symmetry $\TT_{i}=\TT'_{i,1}:\tUB\rightarrow \tUB$ on the generators as follows:
	\begin{align*}
		\TT_i(E_i)=&-(K'_i)^{-1}F_i	,\qquad \TT_i(F_i)=-E_iK_i^{-1},\\
		\TT_i(E_{jl})=&\sum_{r+s=-a_{ij}}(-1)^{r}v^{r}E_i^{(r)}E_{jl}E_i^{(s)}, \text{ for }(j,l)\in\I^{\infty}\text{ and }\; j\neq i,
		\\
		\TT_i(F_{jl})=&\sum_{r+s=-a_{ij}}(-1)^{r}v^{-r}F_i^{(s)}F_{jl}F_i^{(r)}, \text{ for }(j,l)\in\I^{\infty}\text{ and }\; j\neq i,
		\\
		\TT_i(K_\mu)=&K_{\bs_i\mu},\qquad \TT_i(K'_\mu)=K'_{\bs_i\mu},\text{ for any }\mu\in\Z\I.
	\end{align*}
	
	For $i\in \Ire$ and $e \in \{\pm 1\}$,
	we define $\TT'_{i,-1},\TT''_{i,1},\TT''_{i,-1}:\tUB\rightarrow \tUB$  such that
\begin{align}
		\label{eq:TTT}
		\psi \TT_{i,e}' \psi&=\TT_{i,-e}',\qquad \sigma \TT_{i,e}' \sigma = \TT_{i,-e}'',
		\qquad
		\psi \TT''_{i,e} \psi =\TT''_{i,-e}.
	\end{align}

	Similarly to Lemma \ref{lem:Upsilon}, we have the following lemma.
	\begin{lemma}[cf. \text{\cite[Theorem 2.11]{LW22}}]
		\label{lem:UpsilonB}
		Let $\ell$ be a sink of $Q$. Then we have an isomorphism of algebras:
		\begin{align}
			\Gamma_{\ell}:\cs\widetilde{\ch}(\bfk Q)&\stackrel{\cong}{\longrightarrow} \cs\widetilde{\ch}(\bfk \bs_\ell Q)\\\notag
			[M]&\mapsto \sqq^{\langle  \res T_M,\res M\rangle_Q} q^{-\langle T_M,M\rangle}  [F_\ell^+(T_M)]^{-1}\diamond [F_\ell^+(X_M)],
		\end{align}
		where $M\in\cc_{\Z_2}(\rep_\bfk^{\rm nil}(Q))$ and $X_M,T_M\in \cc_{\Z_2}(\ct)$ ($T_M$ acyclic) fit \eqref{resol}.
	\end{lemma}

	\begin{proposition} [\text{\cite[Proposition 2.12]{LW22}}]
		Let $Q$ be a quiver, and $\ell\in Q_0$ be a sink.
		Then the isomorphism $\Gamma_{\ell}:\cs\widetilde{\ch}(\bfk Q)\stackrel{\cong}{\longrightarrow} \cs\widetilde{\ch}(\bfk \bs_\ell Q)$ sends
		\begin{align}
			\Gamma_{\ell}([M])&= [F_{\ell}^+(M)], \quad \forall M\in\cc_{\Z_2}(\ct),
			\label{eqn:reflection 11}
			\\
			\Gamma_{\ell}([C_{S_\ell}])&=
			\sqq[K^*_{S_\ell}]^{-1}* [C^*_{S_{ \ell}}],\qquad \Gamma_{\ell}([C^*_{S_\ell}])=
			\sqq[K_{S_\ell}]^{-1}* [C_{S_{ \ell}}],
			\label{eqn:reflection 21}
			\\
			\Gamma_{\ell}([K_\alpha])&= [K_{\bs_{\ell}\alpha}], \qquad \Gamma_{\ell}([K^*_\alpha])= [K^*_{\bs_{\ell}\alpha}],
			\quad \forall\alpha\in K_0(\rep_\bfk^{\rm nil}(Q)).\label{eqn:reflection 41}
		\end{align}
	\end{proposition}

	Similar to Theorem \ref{thm:main1}, we have the following result, with the help of Lemma \ref{lem:UpsilonB}.
	
	\begin{theorem}
		\label{thm:mainUB}
		Let $Q$ be a quiver, and $\ell\in Q_0$ be a sink. Then we have the following commutative diagram:
		\begin{equation}
			\label{eq:commdiag2}
			\xymatrix{\tUB|_{v=\sqq}\ar[r]^{\TT_\ell} \ar[d]^{\Psi_Q}&  \tUB|_{v=\sqq}\ar[d]^{\Psi_{\bs_\ell Q}} \\
				\cs\widetilde{\ch}(\bfk Q)\ar[r]^{\Gamma_\ell} & \cs\widetilde{\ch}(\bfk \bs_\ell Q)	}
		\end{equation}
	\end{theorem}

	\begin{corollary}
		For any $i\in\Ire$ and $e \in \{\pm 1\}$, we have $\TT'_{i,e},\TT''_{i,e}:\tUB\rightarrow \tUB$ are algebra automorphisms, and $\TT_{i,e}' =(\TT_{i,-e}'')^{-1}$.	
	\end{corollary}
	
	\begin{proof}

		Using Theorem \ref{thm:mainUB} and \eqref{eq:TTT}, we have $\TT'_{i,e},\TT''_{i,e}:\tUB\rightarrow \tUB$ are algebra homomorphisms.
		
		By \eqref{eq:TTT}, we only consider the case $e=1$, i.e., we will show that $\TT_{i,1}' =(\TT_{i,-1}'')^{-1}$.

		Given a quiver $Q$ such that  $i$ is a sink, we have two isomorphisms, $\Gamma_i:\cs\widetilde{\ch}(\bfk Q)\stackrel{\cong}{\longrightarrow} \cs\widetilde{\ch}(\bfk \bs_\ell Q)$ and $\Gamma_i^-:\cs\widetilde{\ch}(\bfk  \bs_\ell  Q)\stackrel{\cong}{\longrightarrow} \cs\widetilde{\ch}(\bfk Q)$, which are inverses to each other; see Lemma \ref{lem:Upsilon} and its proof.
		
		Dually, one can obtain the explicit actions of $\Gamma_i^-$. In particular, similar to Theorem \ref{thm:mainUB}, we have the following commutative diagram:
		\begin{align}
			\label{eq:defT2}
			\xymatrix{
				\tUB|_{v=\sqq}  \ar[r]^{\TT''_{i,-1}}  \ar[d]^{{\Psi}_{\bs_i Q}} & \tUB|_{v=\sqq} \ar[d]^{{\Psi}_{Q}} \\
				\ar[r]^{\Gamma_i^-}  \cs\widetilde{\ch}(\bfk\bs_i Q) & \cs\widetilde{\ch}(\bfk Q)
			}
		\end{align}
		Combining with the above, we conclude that $\TT_{i,1}'$ and $\TT_{i,-1}''$ are inverses to each other.
	\end{proof}
	
	\begin{remark}
		The actions of $\TT'_{i,e}, \TT''_{i,e}$ on $\tU$ factor through the quotient $\UB =\tUB \big/ ( K_i' K_i-1\mid i\in\I )$ to the corresponding automorphisms on $\UB$, and the formulas are then reduced to Lusztig's formulas \cite[\S37.1.3]{Lus93} upon the identification $K_i' =K_i^{-1}$ when $\UB$ is a quantum group.
	\end{remark}

	\begin{proposition}
		The symmetries $\TT'_{i,e},\TT''_{i,e}:\tUB\rightarrow \tUB$ for $i\in\Ire$ satisfy the braid group relations associated to the Borcherds-Cartan datum.
	\end{proposition}
	
	\begin{proof}
		We only need to consider the symmetries $\TT_i=\TT'_{i,1}$.

		The braid group associated to the Borcherds-Cartan datum is generated by $r_i$ for all $i\in\Ire$.
		Denote by $$F(r_i,r_j)=\begin{cases}r_ir_j& \text{ if }a_{ij}=0, \\
			r_ir_jr_i & \text{ if }a_{ij}=-1.\end{cases}$$
		Then the braid group is generated by $r_i$ ($i\in\Ire$) with defining relations
		$F(r_i,r_j)=F(r_j,r_i)$ for any $i\neq j\in\Ire$ if $a_{ij}=0\text{ or }-1$.
		
		Obviously, we have 	
		\begin{align*}F(\TT_i,\TT_j)(K_s)= &F(\TT_j,\TT_i)(K_s),\qquad 	F(\TT_i,\TT_j)(K'_s)= F(\TT_j,\TT_i)(K'_s),\quad \forall s\in\I.
		\end{align*}
		It remains  to check that
		\begin{align}
			\label{eq:braid}
			F(\TT_i,\TT_j)(E_{kl})=& F(\TT_j,\TT_i)(E_{kl}),\qquad F(\TT_i,\TT_j)(F_{kl})= F(\TT_j,\TT_i)(F_{kl}),
		\end{align}
		for any $(k,l)$ with $k\in\I$ and $1\leq l\leq m_k$.
		
		Fix $(k,l)$ with $k\in\I$ and $1\leq l\leq m_k$ in the following.
		Denote by $\tUB_{kl}$ the subalgebra of $\tUB$ generated by
		$$K_s,K_s' \;(s\in\I);\;\; E_{kl}, F_{kl}; \;\;\text{and}\;\; E_{s,1},F_{s,1}\; (k\neq s\in\I).$$
		Then it is obvious that $\TT_i$ preserves $\tUB_{kl}$. Let $\tU'$ be the subalgebra of the quantum Borcherds-Bozec algebra $\tU$ generated by $K_s,K_s'$, $e_{s,1}$, $f_{s,1}$ ($s\in\I$). We also have the braid group action $T_i=T'_{i,1}$ of $\tU$ preserves $\tU'$.
		Then $\tUB_{kl}$ is  isomorphic to $\tU'$ with the isomorphism given by $\phi:\tUB_{kl}\stackrel{\cong}{\longrightarrow}\tU',$
		which sends
		$$K_s\mapsto K_s, K_s'\mapsto K_s'\;(s\in\I); \;\; E_{kl}\mapsto -e_{k,1}, F_{kl}\mapsto \sqq f_{k,1};\;\; E_{s,1}\mapsto -e_{s,1}, F_{s,1}\mapsto \sqq f_{s,1}\; (k\neq s\in\I).$$
		By definition, we have
		$\TT_i:\tUB_{kl}\rightarrow \tUB_{kl}$ coincides with  $T_{i}:\tU'\rightarrow \tU'$ up to the isomorphism $\phi$. So by \cite[Theorem 2.6]{FT21}, we have obtained
		\eqref{eq:braid}.
		
		This proves the proposition.
	\end{proof}

\end{document}